\documentclass[11pt]{amsart}

\usepackage{amscd,amsmath,amssymb,fancyhdr,color}
\usepackage[utf8x]{inputenc}
\usepackage{amsfonts}
\usepackage{amsthm}
\usepackage{accents}
\usepackage{graphicx}
\usepackage{float}
\usepackage{subcaption}
\usepackage{verbatim}
\usepackage{indentfirst}
\usepackage{dsfont}
\usepackage{tikz}
\usetikzlibrary{matrix}
\usepackage{enumerate}
\usepackage{extarrows}

\usepackage{scalerel,stackengine}
\stackMath
\newcommand\reallywidehat[1]{%
	\savestack{\tmpbox}{\stretchto{%
			\scaleto{%
				\scalerel*[\widthof{\ensuremath{#1}}]{\kern-.6pt\bigwedge\kern-.6pt}%
				{\rule[-\textheight/2]{1ex}{\textheight}}
			}{\textheight}%
		}{0.5ex}}%
	\stackon[1pt]{#1}{\tmpbox}%
}

\usepackage{faktor}

\usepackage[backref=page]{hyperref}
\renewcommand*{\backref}[1]{}
\renewcommand*{\backrefalt}[4]{%
    \ifcase #1 (Not cited.)%
    \or        (Cited on page~#2.)%
    \else      (Cited on pages~#2.)%
    \fi}

\hypersetup{
     colorlinks   = true,
     citecolor    = magenta
}

\newcommand{\K}{K\"ahler}

\DeclareMathOperator{\reg}{reg}
\DeclareMathOperator{\supp}{supp}

\numberwithin{equation}{section}

\def\eqref#1{(\ref{#1})}

\newcommand{\del}{\partial}
\newcommand{\delb}{\overline{\partial}}

\def\1{\sqrt{-1}\:}

\newcommand{\cntrct}                
{\hspace{2pt}\raisebox{1pt}{\text{$\lrcorner$}}\hspace{2pt}}

\newcommand{\Aut}{\operatorname{Aut}}

\newcommand{\Sing}{\operatorname{Sing}}

\newcommand{\Deck}{\operatorname{Deck}}

\renewcommand{\Re}{\operatorname{Re}}

\newcommand{\ie}{{\em i.e. }}
\newcommand{\eg}{{\em e.g. }}

\renewcommand{\to}{\longrightarrow}

\newcounter{Mycounter}[section]
\newcounter{lemma}[section]
\newcounter{claim}[section]

\newcounter{sublemma}[section]

\newcounter{corollary}[section]

\newcounter{theorem}[section]

\newcounter{conjecture}[section]

\newcounter{proposition}[section]

\newcounter{definition}[section]

\newcounter{example}[section]

\newcounter{remark}[section]

\newcounter{problem}[section]

\newcounter{question}[section]

\makeatletter

\@addtoreset{equation}{section}

\@addtoreset{footnote}{section}

\makeatother

\usetikzlibrary{arrows,chains,matrix,positioning,scopes}

\makeatletter
\tikzset{join/.code=\tikzset{after node path={%
			\ifx\tikzchainprevious\pgfutil@empty\else(\tikzchainprevious)%
			edge[every join]#1(\tikzchaincurrent)\fi}}}
\makeatother

\tikzset{>=stealth',every on chain/.append style={join},
	every join/.style={->}}

\makeatletter
\newtheorem*{rep@theorem}{\rep@title}
\newcommand{\newreptheorem}[2]{%
	\newenvironment{rep#1}[1]{%
		\def\rep@title{\ref{##1}}%
		\begin{rep@theorem}}%
		{\end{rep@theorem}}}
\makeatother

\newreptheorem{theorem}{Theorem}

\begin{document}
	
\newpage


\title[Coverings of LCK complex spaces]{Coverings of locally conformally K\"ahler complex spaces}

\author{Ovidiu Preda}
\address{Ovidiu Preda \newline
\textsc{\indent Institute of Mathematics ``Simion Stoilow'' of the Romanian Academy\newline 
	\indent 21 Calea Grivitei Street, 010702, Bucharest, Romania}}
\email{ovidiu.preda@imar.ro; ovidiu.preda18@icloud.com}

\author{Miron Stanciu}
\address{Miron Stanciu \newline
\textsc{\indent Institute of Mathematics ``Simion Stoilow'' of the Romanian Academy\newline 
	\indent 21 Calea Grivitei Street, 010702, Bucharest, Romania}}
\email{miron.stanciu@imar.ro; mirostnc@gmail.com}

\thanks{Ovidiu Preda was partially supported by a grant of Ministry of Research and Innovation, CNCS - UEFISCDI, project number
	PN-III-P1-1.1-PD-2016-0182, within PNCDI III. \\\\[.1in]
	{\bf Keywords:} Complex spaces, locally conformally \K .\\
	{\bf 2010 Mathematics Subject Classification:} 32C15, 53C55.
}

\date{\today}

\begin{abstract}
In this paper, we study the properties of coverings of locally conformally \K \ (LCK) spaces with singularities. We begin by proving that a space is LCK if any only if its universal cover is \K, thereby  generalizing a result from \cite{georgeovidiu}. We then show that a complex space which projects over an LCK space with discrete fibers must also carry an LCK structure.
\end{abstract}

\maketitle

\hypersetup{linkcolor=blue}
\tableofcontents

\section{Introduction}

The notion of locally conformally \K \ (LCK) manifolds was first introduced by I. Vaisman \cite{vaisman}. As the name implies, they are manifolds endowed with a complex structure and a non-degenerate two-form which is locally conformal to a \K \ form. They can equivalently be defined as quotients of \K \ manifolds by a discrete group of homotheties of the \K \ form. Since their introduction, they have been extensively studied (see \cite{OV} for an overview of the field and \eg \cite{nico}, \cite{mmo}, \cite{alex} for some recent results).

Our main goal is to study which fundamental properties of LCK structures still hold true in the singular setting. 
The first to introduce the definition of \K \ form on complex spaces with singularities was Grauert in \cite{grauert} and this study was continued by Moishezon in \cite{moish}. Notable results about the stability of \K \ structures were obtained by \cite{var}. Recently, in \cite{georgeovidiu}, the authors gave a natural definition of the same type for LCK spaces and studied the properties of its universal cover, proving, using resolutions of singularities, that, if the space has quotient singularities, the result from the smooth case holds true \ie the universal cover is \K. They conjectured that this should hold without the additional condition. 

One of the key difficulties in working with non-smooth complex spaces is that there is no workable definition for $(p, q)$-forms (for a study of differential forms on singular spaces, see \cite[pp. 375-388]{praga}).

In this paper, we firstly improve the main result from \cite{georgeovidiu} by removing any requirement on the singularity types and in fact finding a way to not mention the singularities at all. Instead, we develop a method for defining closed $1$-forms on complex spaces (without first saying what a $1$-form is) and integrating them along curves. We are then able to essentially use the same proof as in the smooth case for recovering the result we wanted.

Finally, we show that the characterization we found is useful \ie by passing to the universal cover, we can prove new results about coverings of LCK spaces with discrete fibers.

\bigskip

In Section \ref{sec:integrare1forme}, we start by giving the definition of closed $1$-forms on complex spaces and proving that a few basic results about their integration along curves that are true in the smooth case also hold in the singular case.

In Section \ref{sec:acopUnivKahler}, we are then able to improve the result from \cite{georgeovidiu} by completely recovering the characterization of LCK spaces as quotients of \K \ spaces by homotheties of its \K \ metric. More precisely, we prove:

\begin{reptheorem}{thm:LCKacoperiri}
	Let $X$ be a complex space. Then $X$ admits an LCK metric if and only if its universal covering $\tilde{X}$ admits a \K \ metric such that the deck automorphisms act on $\tilde{X}$ by positive homothethies of the \K \ metric.
\end{reptheorem}

In Section \ref{sec:acopVaj}, using the above characterization, we generalize to the locally conformal setting a result by \cite{vaj} about maps with discrete fibers into \K \ spaces. Specifically, we prove:

\begin{reptheorem}{thm:acopDiscerete}
	Let $g: X \to Y$ be a holomorphic map between complex spaces with discrete fibers and assume $(Y, \omega, \theta)$ is LCK. Then $X$ also carries an LCK structure.
\end{reptheorem}

\newpage

\section{Integrating closed 1-forms on complex spaces}
\label{sec:integrare1forme}

Let $X$ be a connected complex space of dimension $n$. 

\begin{definition}
\label{def:1-form}
	\begin{enumerate}[1)]
		\item Denote by 
		\begin{equation*}
		\begin{split}
		\tilde{Z}^1(X) = \{ (U_\alpha, f_\alpha)_{\alpha \in A} \ | \ (U_\alpha)_\alpha \text{ a covering of } X, f_\alpha : U_\alpha \to \mathbb{C} \text { smooth } \\
		\text {and } f_\alpha - f_\beta \text{ locally constant on } U_\alpha \cap U_\beta, \forall \alpha, \beta \in A \}.
		\end{split}
		\end{equation*}
		
		We define an equivalence relation on $\tilde{Z}^1(X)$ by 
		\[
		(U_\alpha, f_\alpha)_{\alpha \in A} \sim (V_\beta, h_\beta)_{\beta \in B} \iff (U_\alpha, f_\alpha)_{\alpha \in A} \cup (V_\beta, h_\beta)_{\beta \in B} \in \tilde{Z}^1(X).
		\]
		
		\item We define the space of \textit{smooth closed 1-forms} on $X$ to be the quotient space $Z^1(X) = \tilde{Z}^1(X)/\sim$. An element $\theta \in Z^1(X)$ is called a \textit{smooth closed 1-form}. 
		
		\item An element $\theta \in Z^1(M)$ is called \textit{exact} if $\theta = \widehat{(X, f)}$ for a smooth $f:X \to \mathbb{C}$. In this case, we make the notation $\theta = df$.
	\end{enumerate}
\end{definition}

\begin{remark}
\label{rem:1-formSing}
	Outside the singular locus $\Sing(X)$, the above definition coincides with the usual one on smooth manifolds.
\end{remark}

Similar to the smooth case, one can define the notion of pullback for closed $1$-forms:

\begin{definition}
\label{def:pullback1forme}
Let $\phi: X \to Y$ be a smooth map between complex spaces and $\theta \in Z^1(Y)$. Let $\theta = \reallywidehat{(U_\alpha, f_\alpha)_{\alpha \in A}}$. 

We denote by $\phi^* \theta = \reallywidehat{(\phi^{-1}(U_\alpha), \phi^*f_\alpha)_{\alpha \in A}} \in Z^1(X)$, and call it \textit{the pullback of} $\theta$ \textit{via} $\phi$. Indeed, one can check that this is well defined and that it coincides with the usual definition on smooth manifolds.
\end{definition}

\begin{proposition}
\label{prop:pullbackComutaCuD}
Let $\phi: X \to Y$ be a smooth map between complex spaces and $f: Y \to \mathbb{C}$ smooth.

Then 
\[
\phi^* df = d \phi^*f. 
\]

\begin{proof}
The statement follows directly from the definitions.
\end{proof}
\end{proposition}

\begin{proposition}
\label{prop:defIntegralei}
	Let $\gamma: [0 ,1] \to X$ be a curve on $X$ and take $\theta \in Z^1(X)$, $\theta = \reallywidehat{(U_\alpha, f_\alpha)_{\alpha \in A}}$.
	
	Choose a partition $0 = a_0 < a_1 < ... < a_m = 1$ such that $\gamma([a_k, a_{k+1}]) \subset U_{\alpha_k}$ for some $\alpha_k \in A$.
	
	Then the number
	\[
	\int_{\gamma, (U_\alpha, f_\alpha)_\alpha, (a_0, ..., a_m)} \theta = \sum_{k=0}^{m} (f_{\alpha_k}(\gamma(a_{k + 1})) - f_{\alpha_k}(\gamma(a_k)))
	\]
	only depends on $\gamma$ and $\theta$.
	
\begin{proof}
	Note that, for any fixed $k$, the choice of $\alpha_k$ such that $\gamma([a_k, a_{k+1}]) \subset U_{\alpha_k}$ does not change the result: if $\gamma([a_k, a_{k+1}]) \subset U_\alpha \cap U_\beta$, then $f_\alpha - f_\beta$ is constant on $\gamma([a_k, a_{k+1}])$, so
	\[
	f_\alpha(\gamma(a_{k+1})) - f_\alpha(\gamma(a_k)) = f_\beta(\gamma(a_{k+1})) - f_\beta(\gamma(a_k)).
	\]
	
	We now prove the independence on the choice of partition of the unit interval. Take two such partitions 
	\begin{equation*}
	\begin{split}
	0 &= a_0 < a_1 < ... < a_m = 1 \\
	0 &= b_0 < b_1 < ... < b_l \ \ = 1
	\end{split}
	\end{equation*}
	and consider the joined partition
	\[
	0 = c_0 < c_1 < ... < c_s = 1.
	\]
	Fix $0 \le k < m$. There exist $0 \le r < p \le s$  such that $c_r = a_k$ and $c_p = a_{k + 1}$. Note that $\gamma([c_t, c_{t + 1}]) \subset U_{\alpha_k}$ for all $r \le t < p$. Then
	\[
	\sum_{t = r}^{p - 1} (f_{\alpha_k} (\gamma(c_{t + 1})) - f_{\alpha_k} (\gamma(c_t))) = f_{\alpha_k}(\gamma(a_{k + 1})) - f_{\alpha_k}(\gamma(a_k))
	\]
	so, according to the remark at the beginning of the proof, and summing by $k$, we have
	\[
	\int_{\gamma, (U_\alpha, f_\alpha)_\alpha, (a_0, ..., a_m)}\theta = \int_{\gamma, (U_\alpha, f_\alpha)_\alpha, (c_0, ..., c_s)}\theta.
	\]
	Similarly, 
	\[
	\int_{\gamma, (U_\alpha, f_\alpha)_\alpha, (b_0, ..., b_l)}\theta = \int_{\gamma, (U_\alpha, f_\alpha)_\alpha, (c_0, ..., c_s)}\theta.
	\]
	
	Lastly, we prove the independence on the choice of representative for $\theta$. Take $\theta = \reallywidehat{(U_\alpha, f_\alpha)_{\alpha \in A}} = \reallywidehat{(V_\beta, h_\beta)_{\beta \in B}}$. Choose a partition $0 = a_0 < ... < a_k = 1$ of the unit interval such that $\gamma([a_k, a_{k + 1}]) \subset U_{\alpha_k} \cap V_{\beta_k}$ for some $\alpha_k \in A, \beta_k \in B$. Then 
	\[
	f_{\alpha_k}(\gamma(a_{k + 1})) - f_{\alpha_k}(\gamma(a_k)) = h_{\beta_k}(\gamma(a_{k + 1})) - h_{\beta_k}(\gamma(a_k)),
	\]
	so 
	\[
	\int_{\gamma, (U_\alpha, f_\alpha)_\alpha, (a_0, ..., a_m)}\theta = \int_{\gamma, (V_\beta, h_\beta)_\beta, (a_0, ..., a_m)}\theta.
	\]
\end{proof}
\end{proposition}

\begin{definition}
\label{def:integrala}
For any curve $\gamma$ and any $1$-form $\theta$, we call the number defined \textit{via} \ref{prop:defIntegralei} $\textit{the integral of } \theta \textit{ along } \gamma$ and denote it by
\[
\int_\gamma \theta.
\]
\end{definition}

\begin{remark}
\label{rem:intSing}
	If the image of $\gamma$ does not intersect the singular locus, the above definition coincides with the regular integral of a $1$-form along a curve on a smooth manifold.
\end{remark}

\medskip

\ref{prop:defIntegralei} has the following immediate

\begin{corollary}
\label{cor:intPuncte}
If $\gamma_0, \gamma_1$ are two curves with the same endpoints in $X$, $\theta \in Z^1(X)$ and there exist $\theta = \reallywidehat{(U_\alpha, f_\alpha)_{\alpha \in A}}$ and a partition $0 = a_0 <a_1 < ... < a_m = 1$ with $\gamma_0([a_k, a_{k+1}]) \ \cup \ \gamma_1([a_k, a_{k+1}]) \subset U_{\alpha_k}$ for some $\alpha_k \in A$ such that $\gamma_0(a_k) = \gamma_1(a_k)$ for all $0 \le k \le m$, then
\[
\int_{\gamma_0} \theta = \int_{\gamma_1} \theta.
\]
	
\end{corollary}

\begin{proposition}
\label{prop:invarOmotopie}
The integral of $\theta$ along $\gamma$ only depends on the homotopy class (with fixed endpoints) of $\gamma$.
\begin{proof}
	Consider a homotopy $\gamma_\cdot: [0,1] \times [0,1] \to X, \gamma_t(0) = a, \gamma_t(1) = b$. We show that 
	\[
	\left( [0, 1] \ni t \to \int_{\gamma_t} \theta \right) \text{ is locally constant.}
	\]
	Fix $t_0 \in [0,1]$. Take $\theta = \reallywidehat{(U_\alpha, f_\alpha)_{\alpha \in A}}$ with $U_\alpha$ coordinate charts and choose a partition $0 = a_0 < ... < a_m = 1$ and an $\epsilon > 0$ such that, for all $0 \le k \le m - 1$,
	\[
	\gamma([a_k, a_{k + 1}] \times (t_0 - \epsilon, t_0 + \epsilon) ) \subset U_{\alpha_k} \text { for some } \alpha_k \in A.
	\]
	
	By \ref{prop:defIntegralei}, we have:
	\begin{multline*}	
	\int_{\gamma_{t}}\theta - \int_{\gamma_{t_0}}\theta = \hspace*{\fill} \\
	\hspace*{\fill} \sum_{k=0}^{m-1}\left(f_{\alpha_k}(\gamma_t(a_{k+1}))-f_{\alpha_k}(\gamma_t(a_k))\right) - \sum_{k=0}^{m-1}\left(f_{\alpha_k}(\gamma_{t_0}(a_{k+1}))-f_{\alpha_k}(\gamma_{t_0}(a_k))\right) =  \\
	\hspace*{\fill} \sum_{k=0}^{m-1}\left( f_{\alpha_k}(\gamma_t(a_{k+1}))-f_{\alpha_k}(\gamma_{t_0}(a_{k+1}))  \right) -
	\sum_{k=0}^{m-1}\left( f_{\alpha_k}(\gamma_t(a_k))-f_{\alpha_k}(\gamma_{t_0}(a_k))  \right) = \\
	\hspace*{\fill} \sum_{k=0}^{m-1}\left( f_{\alpha_{k+1}}(\gamma_t(a_{k+1}))-f_{\alpha_{k+1}}(\gamma_{t_0}(a_{k+1}))  \right) -
	\sum_{k=0}^{m-1}\left( f_{\alpha_k}(\gamma_t(a_k))-f_{\alpha_k}(\gamma_{t_0}(a_k))  \right)=\\ 
	=0, \mbox{\ for any\ } t\in(t_0-\epsilon,t_0+\epsilon). \hspace*{\fill}
	\end{multline*}
	
	For the third equality, replacing $f_{\alpha_k}$ with $f_{\alpha_{k+1}}$ in the first sum is possible, since the curve $s\mapsto \gamma_s(a_{k+1})$, with $s$ between $t$ and $t_0$, is contained entirely in $U_{\alpha_k} \cap U_{\alpha_{k+1}}$.
\end{proof}
\end{proposition}

\medskip

\begin{corollary}
\label{cor:Poincare}
	If $X$ is simply connected, every closed $1$-form $\theta$ is exact.
	
\begin{proof}
	Take $\theta = \reallywidehat{(U_\alpha, f_\alpha)_{\alpha \in A}}$. Fix $x_0 \in X$. We define $f: X \to \mathbb{C}$,
	\[
	f(x) = \int_{x_0}^x \theta := \int_\gamma \theta, \text{ where } \gamma \text{ is a curve from } x_0 \text{ to } x.
	\]
	Since $X$ is simply connected and by \ref{prop:invarOmotopie}, this does not depend on the choice of $\gamma$.
	
	We now show that $\theta = df$ \ie $f - f_\alpha$ is locally constant on $U_\alpha$ for all $\alpha \in A$. Choose $x_\alpha \in U_\alpha$. Then
	\[
	f(x) = \int_{x_0}^{x_\alpha} \theta + \int_{x_\alpha}^{x} \theta = \int_{x_0}^{x_\alpha} \theta + f_\alpha(x) - f_\alpha(x_\alpha),
	\]
	so
	\[
	f(x) - f_\alpha(x) = \int_{x_0}^{x_\alpha} \theta - f_\alpha(x_\alpha), \ \forall x \in U_\alpha,
	\]
	and the right-hand side is a constant which does not depend on the choice of $x_\alpha \in U_\alpha$.
\end{proof}
\end{corollary}

\section{Locally conformally \K \ spaces}
\label{sec:acopUnivKahler}

In this section, we give the definitions of the main objects we will be working with, as well as prove a result characterizing LCK spaces \textit{via} their universal covering, thereby generalizing a well-known property which is true in the smooth case.

 The following definition is very slightly adapted from \cite{moish} to resemble \ref{def:1-form}:

\begin{definition}
\label{def:Kahler}
	\begin{enumerate}[1)]
		\item Denote by 
		\begin{equation*}
		\begin{split}
		\tilde{\mathcal{K}}(X) = \{ & (U_\alpha, \varphi_\alpha)_{\alpha \in A} \ | \ (U_\alpha)_\alpha \text{ a covering of } X, \varphi_\alpha : U_\alpha \to \mathbb{R} \\ & \text { strongly plurisubharmonic and} \ \varphi_\alpha - \varphi_\beta = \Re g_{\alpha \beta} \text { on } U_\alpha \cap U_\beta, \\ & \text{ for some holomorphic } g_{\alpha \beta}: U_\alpha \cap U_\beta \to \mathbb{C}, \forall \alpha, \beta \in A \}.
		\end{split}
		\end{equation*}
		
		We define an equivalence relation on $\tilde{\mathcal{K}}(X)$ by 
		\[
		(U_\alpha, \varphi_\alpha)_{\alpha \in A} \sim (V_\beta, \psi_\beta)_{\beta \in B} \iff (U_\alpha, \varphi_\alpha)_{\alpha \in A} \cup (V_\beta, \psi_\beta)_{\beta \in B} \in \tilde{K}(X).
		\]
		
		\item We define the space of \textit{K\"ahler forms} on $X$ to be the quotient space $\mathcal{K}(X) = \tilde{\mathcal{K}}(X)/\sim$. An element $\omega \in \mathcal{K}(X)$ is called a \textit{K\"ahler form}.
		
		\item $(X, \omega)$ is called a \textit{K\"ahler space}. 
	\end{enumerate}	
	
\end{definition}

\begin{remark}
\label{rem:KahlerPeReg}
One can also give an equivalent definition for a \K \ form on $X$ \textit{via} differential conditions on the regular locus: a \K \ form is given by $(U_\alpha, \varphi_\alpha)_\alpha$ with $\varphi_\alpha$ strongly plurisubharmonic on $U_\alpha$ if any only if 
\[
i \del \delb (\varphi_\alpha - \varphi_\beta) = 0
\]
on $U_\alpha \cap U_\beta \cap X_{\reg}$.

\end{remark}

\bigskip

Inspired by the characterization given in \ref*{rem:KahlerPeReg}, the following definition for LCK spaces was first introduced in \cite{georgeovidiu}, adapting the most well-suited of the equivalent definitions of LCK manifolds. For technical reasons, we introduce at the same time what we call \textit{locally conformally pre\K} metrics.

\newpage

\begin{definition}
	\label{def:LCK}
	\begin{enumerate}[1)]
		\item Denote by 
		\begin{equation*}
		\begin{split}
		\widetilde{\mathcal{LCK}}(X) = \{ & (U_\alpha, \varphi_\alpha, f_\alpha)_{\alpha \in A} \ | \ (U_\alpha)_\alpha \text{ a covering of } X, \varphi_\alpha : U_\alpha \to \mathbb{R} \\ & \text { strongly plurisubharmonic, } f_\alpha: U_\alpha \to \mathbb{R} \text{ smooth and} \\ & 	e^{f_\alpha} i \del \delb \varphi_\alpha = e^{f_\beta} i \del \delb \varphi_\beta \text{ on } U_\alpha \cap U_\beta \cap X_{\reg}, \forall \alpha, \beta \in A \}.
		\end{split}
		\end{equation*}
		
		\begin{equation*}
		\begin{split}
		\widetilde{\mathcal{LC}p\mathcal{K}}(X) = \{ & (U_\alpha, \varphi_\alpha, f_\alpha)_{\alpha \in A} \ | \ (U_\alpha)_\alpha \text{ a covering of } X, \varphi_\alpha : U_\alpha \to \mathbb{R} \\ & \text { plurisubharmonic, } f_\alpha: U_\alpha \to \mathbb{R} \text{ smooth and} \\ & 	e^{f_\alpha} i \del \delb \varphi_\alpha = e^{f_\beta} i \del \delb \varphi_\beta \text{ on } U_\alpha \cap U_\beta \cap X_{\reg}, \forall \alpha, \beta \in A \}.
		\end{split}
		\end{equation*}
		
		We define equivalence relations on $\widetilde{\mathcal{LCK}}(X)$ and $\widetilde{\mathcal{LC}p\mathcal{K}}(X)$ by 
		\begin{equation*}
		\begin{split}
		(U_\alpha, \varphi_\alpha, f_\alpha)_{\alpha \in A} \sim (V_\beta, \psi_\beta, h_\beta)_{\beta \in B} \iff \\ (U_\alpha, \varphi_\alpha, f_\alpha)_{\alpha \in A} \cup (V_\beta, \psi_\beta, h_\beta)_{\beta \in B} \in \widetilde{\mathcal{LCK}}(X).
		\end{split}
		\end{equation*}
		and similarly,
		\begin{equation*}
		\begin{split}
		(U_\alpha, \varphi_\alpha, f_\alpha)_{\alpha \in A} \sim (V_\beta, \psi_\beta, h_\beta)_{\beta \in B} \iff \\ (U_\alpha, \varphi_\alpha, f_\alpha)_{\alpha \in A} \cup (V_\beta, \psi_\beta, h_\beta)_{\beta \in B} \in \widetilde{\mathcal{LC}p\mathcal{K}}(X).
		\end{split}
		\end{equation*}
		
		\item We define the space of \textit{locally conformally K\"ahler (LCK) forms} (or \textit{metrics}) on $X$ to be the quotient space $\mathcal{LCK}(X) = \widetilde{\mathcal{LCK}}(X)/\sim$ and the space of \textit{locally conformally preK\"ahler (LCpK) forms} (or \textit{metrics}) on $X$ to be the quotient space $\mathcal{LC}p\mathcal{K}(X) = \widetilde{\mathcal{LC}p\mathcal{K}}(X)/\sim$. An element $\omega \in \mathcal{LC}(p)\mathcal{K}(X)$ is called a LC(p)K metric. Obviously, $\mathcal{LCK}(X) \subset \mathcal{LC}p\mathcal{K}(X)$.
		
		\item The covering $(U_\alpha)_\alpha$ together with the functions $f_\alpha$ give rise to a closed $1$-form $\theta$ on $X$, as described in \ref{def:1-form}. We call $\theta$ the \textit{Lee form} associated to the LC(p)K metric.
		
		\item $(X, \omega, \theta)$ is called a \textit{LC(p)K space}.
		
		\item An LCpK metric $\omega$ is called \textit{with potential} if $\omega = \reallywidehat{(X, \varphi, f)}$ for $\varphi:X \to \mathbb{\mathbb{R}}$ plurisubharmonic and $f:X \to \mathbb{R}$ smooth. In this case, we make the notation $\omega = e^f \del \delb \varphi$.
	\end{enumerate}	
	
\end{definition}

\begin{remark}
\label{rem:K&LCKpeReg}
	On $X_{\reg}$, \ref{def:Kahler} and \ref{def:LCK} revert to the classical definitions for \K, locally conformally \K \ and locally conformally pre\K \ manifolds, respectively.
\end{remark}

\begin{remark}
\label{rem:LCKinvarLaTransConf}
	Let $(X, \omega, \theta)$ be an LC(p)K space and $f: X \to \mathbb{R}$ a smooth map. Let $\omega$ be given by $(U_\alpha, \varphi_\alpha, f_\alpha)$. Then $X$ is also LC(p)K with the metric given by $(U_\alpha, \varphi_\alpha, f + f_\alpha)$, and we denote it by $(X, e^f \omega, \theta + df)$. Indeed, one can check that its Lee form is $\theta + df$, with the notations of \ref{def:1-form}.
\end{remark}

\begin{remark}
\label{rem:GCK}
As the notations suggest, if $(X, \omega, \theta)$ is LCK and $\theta = df$ \ie is exact, then $(X, e^{-f}\omega, \theta - df = 0)$ is \K. In this case, we call $(X, \omega, \theta)$ \textit{globally conformally K\"ahler (GCK)}. 
\end{remark}

\medskip

\begin{definition}
	\label{def:pullbackaLCK}
	Let $\phi: X \to Y$ be a holomorphic map between complex spaces and assume $(Y, \omega, \theta)$ is LCpK. Let $\omega = \reallywidehat{(U_\alpha, \varphi_\alpha, f_\alpha)_{\alpha \in A}}$. 
	
	Then one sees immediately that $(\phi^{-1}(U_\alpha), \phi^* \varphi_\alpha, \phi^*f_\alpha)_{\alpha \in A}$ gives an LCpK metric on $X$, which we denote by $\phi^* \omega$, with Lee form $\phi^* \theta$. We call $\phi^* \omega$ the \textit{pullback of $\omega$ via $\phi$}.
\end{definition}

\medskip

As opposed to LCpK forms or closed $1$-forms, the pullback of LCK metrics along a holomorphic map is not necessarily LCK. However, there is a particular case where this is true:

\begin{proposition}
\label{prop:pullbackLCK}
Let $\phi: X \to Y$ be a a local biholomorphism between complex spaces and assume $(Y, \omega, \theta)$ is LCK. Then $\phi^* \omega$ is LCK, with Lee form $\phi^* \theta$.

\begin{proof}
The result is immediate. As $\phi$ is a local biholomorphism, $\phi^* \varphi_\alpha$ are also strictly plurisubharmonic. In the same way, if $g_{\alpha \beta}$ is holomorphic on $U_\alpha \cap U_\beta$, then $\phi^* g_{\alpha \beta}$ is holomorphic on $\phi^{-1}(U_\alpha) \cap \phi^{-1}(U_\beta)$.
\end{proof}
\end{proposition}

\begin{definition}
	Let $(X, \omega, \theta)$ be an LC(p)K space. An automorphism $\gamma \in \Aut(X)$ acts by homothethies of $\omega$ if $\gamma^* \omega = e^c \omega$ for some $c \in \mathbb{R}$.
\end{definition}

\bigskip

We now prove a result giving an alternative characterization for LCK spaces which is a natural extension of the one true in the smooth case and has virtually the same proof, with the help of the notions developed in Section \ref{sec:integrare1forme}. This generalizes a result in \cite{georgeovidiu}:

\begin{theorem}
\label{thm:LCKacoperiri} 
	Let $X$ be a complex space. Then $X$ admits an LCK metric if and only if its universal covering $\tilde{X}$ admits a \K \ metric such that the deck automorphisms act on $\tilde{X}$ by positive homothethies of the \K \ metric.

\begin{proof}
	Assume that $(X, \omega, \theta)$ is an LCK space and denote by $\pi : \tilde{X} \to X$ the canonical projection. Then $(\tilde{X}, \pi^* \omega, \pi^* \theta)$ is an LCK space by \ref{prop:pullbackLCK}. Futhermore, by \ref{cor:Poincare}, $\theta = df$ for some smooth $f: \tilde{X} \to \mathbb{R}$, so \ref{rem:GCK} implies that $(\tilde{X}, e^{-f} \pi^* \omega)$ is \K. 
	
	Lastly, for any $\gamma$ a deck automorphism of $\tilde{X} \to X$, we have $\theta = \gamma^* \theta$, so, according to \ref{prop:pullbackComutaCuD}, $df = \gamma^* df = d \gamma^* f$, so $\gamma^* f - f$ is constant. Hence, if 
	
	Conversely, assume $(\tilde{X}, \tilde{\omega})$ is \K, with the deck group acting by positive homothethies of $\tilde{\omega}$ on it. Then the map
	\[
	\chi: \Deck(\tilde{X}/X) = \pi_1(X) \to \mathbb{R}_{>0}, \ \chi(\gamma) = \frac{\gamma^* \tilde{\omega}}{\tilde{\omega}}
	\]
	is a character of $\pi_1(X)$. We can thus construct the line bundle $L = \tilde{X} \times_\chi \mathbb{R}$ over $X$ \ie
	\[
	L = (\tilde{X} \times \mathbb{R}) / ((\tilde{x}, a) \sim (\gamma(\tilde{x}), \chi(\gamma) a)).
	\]
	Since there exists a covering of $X$ and a choice of transition maps all of which are positive, $L$ is trivial. Consider a section $v$ of $L$ which is nonzero at every point. 	
	Take $\tilde{v}=\pi^\star(v)$. Since $L$ is trivial, $\pi^* L \to \tilde{X}$ is also trivial, hence $\tilde{v}$ can be viewed as a function $\tilde{v}:\tilde{X} \rightarrow \mathbb{R}_{>0}$, which can be written $\tilde{v}=e^{-\tilde{f}}$.
	By the definition of $L$, we have: 
	\[
	\frac{\gamma^* \tilde{v}}{\tilde{v}}=\chi(\gamma) = \frac{\gamma^* \tilde{\omega}}{\tilde{\omega}}
	\]
	Consequently, $e^{\tilde{f}}\tilde{\omega}$ is deck invariant, and it descends to a metric $\omega$ on $X$, which, by its definition, is LCK.
\end{proof}

\end{theorem}

\section{Mappings with discrete fibers between LCK spaces}
\label{sec:acopVaj}

In this section, we prove the following 

\begin{theorem}
	\label{thm:acopDiscerete}
	Let $g: X \to Y$ be a holomorphic map between complex spaces with discrete fibers and assume $(Y, \omega, \theta)$ is LCK. Then $X$ also carries an LCK structure. 
\end{theorem}

This generalizes a result for \K \ spaces, see \cite[Theorem 2]{vaj}. In fact, our proof consists of a careful intermix of the proof of that result with further topological arguments.

\medskip

We begin by stating an easy topological fact:

\begin{lemma}
	\label{lem:acopTop}
	Let $g:X \to Y$ be a continuous function with discrete fibers between locally compact topological spaces and choose $x \in U \subset X$ with $U$ open. 
	
	Then there exists an open set $W \ni g(x)$ such that the connected component of $g^{-1}(W)$ containing $x$ is contained in $U$.
\end{lemma}

\begin{definition}
	\label{def:acopVaj}
	Let $g: X \to (Y, \omega, \theta)$ be a holomorphic map between complex spaces with discrete fibers, where $(Y, \omega, \theta)$ is LCK.
	
	We call the coverings $(U_\alpha)_{\alpha \in A}$ of $X$ and $(V_\alpha)_{\alpha \in A}$ of $Y$ \textit{well-related} if:
	\begin{enumerate}[1)]
		\item There exist $\psi_\alpha : V_\alpha \to \mathbb{R}$ strongly plurisubharmonic and $f_\alpha: V_\alpha \to \mathbb{R}$ smooth, such that $\omega = \reallywidehat{(V_\alpha, \varphi_\alpha, f_\alpha)_{\alpha \in A}}$.
		
		\item The covering $(U_\alpha)_{\alpha \in A}$ is locally finite and relatively compact.
		
		\item There exist functions $\varphi_\alpha: U_\alpha \to (0, \infty)$ which are strongly plurisubharmonic.
		
		\item Each $U_\alpha$ is a connected component of $g^{-1}(V_\alpha)$.
	\end{enumerate}
\end{definition}

The existence of such well-related coverings follows from \ref{lem:acopTop}.

\medskip

As we want to use the \K ianity of the universal covers of $X$ and $Y$ given by \ref{thm:LCKacoperiri}, we now prove the following technical

\begin{proposition}
\label{prop:acopeririVajUniversala}
	Let $g: X \to (Y, \omega, \theta)$ be a holomorphic map between complex spaces with discrete fibers, where $(Y, \omega, \theta)$ is LCK, and $(U_\alpha)_{\alpha \in A}$ and $(V_\alpha)_{\alpha \in A}$ well-related coverings of $X$ and $Y$ respectively.
	
	Denote by $\pi_X : \tilde{X} \to X$ and $\pi_Y : \tilde{Y} \to Y$ the universal coverings of $X$ and $Y$ and consider $\tilde{g}:\tilde{X} \to \tilde{Y}$ the lifting of $g$, as in the following commutative diagram:
	\begin{figure}[H]
		\centering
		\label{diag:thmEmbeddingOutline}
		\begin{tikzpicture}
		\matrix (m) [matrix of math nodes,row sep=3em,column sep=4em,minimum width=2em]
		{
			\tilde{X} & \tilde{Y} \\
			X &  (Y, \omega, \theta) \\};
		\path[-stealth]
		(m-1-1) edge node [above] {$\tilde{g}$} (m-1-2)
		(m-1-1) edge node [left] {$\pi_X$} (m-2-1)
		(m-1-2) edge node [right] {$\pi_Y$} (m-2-2)
		(m-2-1) edge node [above] {$g$} (m-2-2);
		\end{tikzpicture}
	\end{figure}
	Notice that $\tilde{g}$ also has discrete fibers.
	
	We may also assume that $U_\alpha$ are small enough such that $\pi_X^{-1}(U_\alpha)$ is a disjoint union of copies of $U_\alpha$. For each $\alpha \in A$, choose an isomorphism between the fiber above elements of $U_\alpha$ and $\Gamma = \Deck(\tilde{X}/X)$, such that 
	\[
	\pi_X^{-1}(U_\alpha) = \coprod\limits_{\gamma \in \Gamma} \tilde{U}_{(\alpha, \gamma)},
	\]
	with $\pi_{X| \tilde{U}_{(\alpha, \gamma)}}: \tilde{U}_{(\alpha, \gamma)} \to U_\alpha$ a biholomorphism. 
	
	For $(\alpha, \gamma) \in A \times \Gamma$, define $\tilde{V}_{(\alpha, \gamma)} = \pi_Y^{-1}(V_\alpha)$.
	
	Then the coverings $(\tilde{U}_{(\alpha, \gamma)})_{(\alpha, \gamma) \in A \times \Gamma}$ of $\tilde{X}$ and $(\tilde{V}_{(\alpha, \gamma)})_{(\alpha, \gamma) \in A \times \Gamma}$ of $\tilde{Y}$ are well-related, where $(\tilde{Y}, e^{-f} \pi_Y^* \omega)$ is \K, with $df = \pi_Y^* \theta$, as proven by \ref{thm:LCKacoperiri}.
	
	Moreover, given a representation $\rho: \Gamma \to \mathbb{R}_+$, we may choose $\tilde{\varphi}_{(\alpha, \gamma)}: \tilde{U}_{(\alpha, \gamma)} \to (0, \infty)$ strongly plurisubharmonic such that 
	\[
	\eta^* \tilde{\varphi}_{(\alpha, \gamma)} = \rho(\eta) \tilde{\varphi}_{(\alpha, \eta^{-1} \gamma)}, \ \forall (\alpha, \gamma) \in A \times \Gamma, \ \forall \eta \in \Gamma.
	\]
	
	\begin{proof}
	We check the conditions imposed by \ref{def:acopVaj}:
		
	\begin{enumerate}[1)]
		\item As shown in \ref{prop:pullbackLCK} and \ref{thm:LCKacoperiri}, the choices
		\[
		\tilde{\psi}_{(\alpha, \gamma)} = \pi_Y^* \psi_\alpha: \tilde{V}_{(\alpha, \gamma)} \to \mathbb{R}
		\]
		are strongly plurisubharmonic and give rise to the \K \ metric $e^{-f} \pi_Y^* \omega$ on $\tilde{Y}$.
			
		\item The covering $(\tilde{U}_{(\alpha, \gamma)})_{(\alpha, \gamma) \in A \times \Gamma}$ is immediately seen to be locally finite and relatively compact.
		
		\item By the definition of well-relatedness, there exist functions $\varphi_\alpha: U_\alpha \to (0, \infty)$ strongly plurisubharmonic. Define
		\[
		\tilde{\varphi}_{(\alpha, \gamma)}: \tilde{U}_{(\alpha, \gamma)} \to (0, \infty), \ \tilde{\varphi}_{(\alpha,  \gamma)} = \rho(\gamma) \pi_X^* \varphi_\alpha.
		\]
		These functions are obviously strongly plurisubharmonic. Moreover, for an $\eta \in \Gamma$, $\eta^* \tilde{\varphi}_{(\alpha, \gamma)} \in \mathcal{C}^\infty(\tilde{U}_{(\alpha, \eta^{-1}\gamma)})$
		\begin{equation*}
		\begin{split}
		\eta^* \tilde{\varphi}_{(\alpha, \gamma)} &= \rho(\gamma) \eta^* \pi_X^* \varphi_\alpha = \rho(\gamma) \pi_X^* \varphi_\alpha = \rho(\eta) \rho(\eta^{-1} \gamma) \pi_X^* \varphi_\alpha \\ &= \rho(\eta) \tilde{\varphi}_{(\alpha, \eta^{-1}\gamma)}.
		\end{split}
		\end{equation*}
		
		\item We need to show that $\tilde{U}_{(\alpha, \gamma)}$ is a connected component of $\tilde{g}^{-1}(\tilde{V}_{(\alpha, \gamma)})$, for each $(\alpha, \gamma) \in A \times \Gamma$. Indeed, $\tilde{g}^{-1}(\tilde{V}_{(\alpha, \gamma)}) = \pi_X^{-1} (g^{-1}(V_\alpha))$ and, by definition, $U_\alpha$ is a connected component of $g^{-1}(V_\alpha)$.
	\end{enumerate}
	\end{proof}
\end{proposition}

\bigskip

We can now turn to the 

\begin{proof}[\textbf{Proof of \ref{thm:acopDiscerete}}]
	Denote again by $\pi_X : \tilde{X} \to X$ and $\pi_Y : \tilde{Y} \to Y$ the universal coverings of $X$ and $Y$ and consider $\tilde{g}:\tilde{X} \to \tilde{Y}$ the lifting of $g$. By \ref{thm:LCKacoperiri}, $(\tilde{Y}, e^{-f} \pi_Y^*\omega)$ is \K, where $df = \pi_Y^*\theta$. 
	
	Pick $\eta \in \Gamma = \Deck(\tilde{X}/X)$. Then
	\begin{equation}
	\label{eq:thm:acopDiscrete1}
	\eta^* \tilde{g}^* (e^{-f}\pi_Y^* \omega) = e^{-\eta^* \tilde{g}^* f} \eta^* \tilde{g}^* \pi_Y^* \omega
	\end{equation}
	(note that this is an equality of LCpK metrics).
	
	But 
	\begin{equation}
	\label{eq:thm:acopDiscrete2}
	\eta^* \tilde{g}^* \pi_Y^* \omega = (\pi_Y \tilde{g} \eta)^* \omega = (g \pi_X \eta)^* \omega = (g \pi_X)^* \omega = (\pi_Y \tilde{g})^* \omega = \tilde{g}^* \pi_Y^* \omega.
	\end{equation}
	
	Similarly, 
	\begin{equation*}
	d(\eta^* \tilde{g}^* f) = \eta^* \tilde{g}^* df = \eta^* \tilde{g}^* \pi_Y^* \theta = \tilde{g}^* \pi_Y^* \theta = \tilde{g}^* df = d(\tilde{g}^* f),
	\end{equation*}
	so $\eta^* \tilde{g}^* f = \tilde{g}^* f + c_\eta$ and thus
	\begin{equation}
	\label{eq:thm:acopDiscrete3}
	e^{-\eta^* \tilde{g}^* f} = e^{-c_\eta} e^{-\tilde{g}^* f}.
	\end{equation}
	Using \eqref{eq:thm:acopDiscrete2} and \eqref{eq:thm:acopDiscrete3} in \eqref{eq:thm:acopDiscrete1}, we conclude that
	\begin{equation}
	\label{eq:thm:acopDiscrete4}
	\eta^* \tilde{g}^* (e^{-f}\pi_Y^* \omega) = e^{-c_\eta} \tilde{g}^* (e^{-f}\pi_Y^* \omega).
	\end{equation} 
	Note also that $\eta \mapsto e^{-c_\eta}$ is a representation of $\Gamma$, which we denote by $\rho$.
	
	\medskip
	
	For $g: X \to (Y, \omega, \theta)$, pick $(U_\alpha)_{\alpha \in A}$ and $(V_\alpha)_{\alpha \in A}$ well-related coverings of $X$ and $Y$ respectively. For $\rho: \Gamma \to \mathbb{R}_+$ chosen as above, we use \ref{prop:acopeririVajUniversala} to construct the well-related coverings $(\tilde{U}_{(\alpha, \gamma)})_{(\alpha, \gamma) \in A \times \Gamma}$ of $\tilde{X}$ and $(\tilde{V}_{(\alpha, \gamma)})_{(\alpha, \gamma) \in A \times \Gamma}$ of $\tilde{Y}$.
	
	The proof now follows the same broad steps as in \cite{vaj}. Consider $K_\alpha \subset U_\alpha$ compact sets with $\bigcup\limits_{\alpha \in A} K_\alpha = X$ and choose $\tau_\alpha \in \mathcal{C}_0^\infty(V_\alpha), \tau_\alpha \ge 0$ and $\tau_\alpha$ equals $1$ on a neighborhood of $g(K_\alpha)$. Let
	\[
	\chi_{(\alpha, \gamma)} = 
	\left\{
	\begin{array}{ll}
	(\tau_\alpha \circ g\circ \pi_X)^2=(\tau_\alpha \circ \pi_Y \circ \tilde{g})^2 & \text{on} \ U_{(\alpha, \gamma)}, \\
	0 & \text{on} \ \tilde{X} \setminus U_{(\alpha, \gamma)}.
	\end{array}
	\right.
	\]
	Note that $\chi_{(\alpha, \gamma)} \varphi_{(\alpha, \gamma)}$ is smooth on $\tilde{X}$ with compact support in $U_{(\alpha, \gamma)}$ and that, for any $\eta \in \Gamma$, $\eta^* \chi_{(\alpha, \gamma)} = \chi_{(\alpha, \eta^{-1}\gamma)}.$ 
	
	Then for any $\epsilon = (\epsilon_\alpha)_{\alpha \in A}$ positive numbers, we define the smooth function $\varphi_\epsilon: \tilde{X} \to (0, \infty)$ by
	\begin{equation*}
	\varphi_\epsilon = \sum\limits_{\substack{\alpha \in A \\ \gamma \in \Gamma}} \epsilon_\alpha \chi_{(\alpha, \gamma)} \varphi_{(\alpha, \gamma)}.
	\end{equation*} 
	For an $\eta \in \Gamma$, we have, by \ref{prop:acopeririVajUniversala},
	\begin{equation}
	\label{eq:thm:acopDiscrete5}
	\begin{split}
	\eta^* \varphi_\epsilon &= \sum\limits_{\substack{\alpha \in A \\ \gamma \in \Gamma}} \epsilon_\alpha \eta^* \left( \chi_{(\alpha, \gamma)} \varphi_{(\alpha, \gamma)} \right) = \rho(\eta) \sum\limits_{\substack{\alpha \in A \\ \gamma \in \Gamma}} \epsilon_\alpha \chi_{(\alpha, \eta^{-1}\gamma)} \varphi_{(\alpha, \eta^{-1}\gamma)} \\ &= \rho(\eta) \varphi_\epsilon.
	\end{split}
	\end{equation} 	
	
	Now consider 
	\[
	\phi_{(\alpha, \gamma)}^\epsilon : \tilde{U}_{(\alpha, \gamma)} \to \mathbb{R}, \ \phi_{(\alpha, \gamma)}^\epsilon = \tilde{g}^* \psi_{(\alpha, \gamma)} + \varphi_\epsilon.
	\]
	Since $\phi_{(\alpha, \gamma)}^\epsilon - \phi_{(\beta, \sigma)}^\epsilon = \tilde{g}^* (\psi_{(\alpha, \gamma)} - \psi_{(\beta, \sigma)})$, these differences are pluriharmonic. 
	
	The fact that one can choose $\epsilon$ and an open covering $(U_{(\alpha, \gamma)}^\prime)_{(\alpha, \gamma)\in A\times \Gamma}$ for which $U_{(\alpha, \gamma)}^\prime \subset U_{(\alpha, \gamma)}$ such that $(U_{(\alpha, \gamma)}^\prime, \phi_{(\alpha, \gamma)}^\epsilon)_{(\alpha, \gamma) \in A \times \Gamma}$ gives a \K \ form on $\tilde{X}$ now follows from the following Claim, which for clarity we will prove at the end:
	
	\medskip
	
	\textbf{Claim:} There exists a set $\epsilon^0 = (\epsilon_\alpha^0)_{\alpha \in A}$ of positive constants such that $\phi_{(\alpha, \gamma)}^\epsilon$ is strongly plurisubharmonic on the support of $\chi_{(\alpha, \gamma)}$ for all $(\alpha, \gamma) \in A \times \Gamma$ and for all $0 < \epsilon < \epsilon^0$.
	
	\medskip
	
	Choose an appropriate $\epsilon$ and denote by $\tilde{\omega}_X$ = $\reallywidehat{(U_{(\alpha, \gamma)}^\prime, \phi_{(\alpha, \gamma)}^\epsilon)_{(\alpha, \gamma) \in A \times \Gamma}}$ the \K \ form on $\tilde{X}$. 
	
	For any $\eta \in \Gamma$, we have, by \eqref{eq:thm:acopDiscrete4} and \eqref{eq:thm:acopDiscrete5}, 
	\[
	\eta^* \tilde{\omega}_X = \eta^* (\tilde{g}^* e^{-f} \pi_Y^* \omega + \del \delb \varphi_\epsilon) = \rho(\eta) (\tilde{g}^* e^{-f} \pi_Y^* \omega + \del \delb \varphi_\epsilon) = \rho(\eta) \tilde{\omega}_X,
	\]
	so $\Gamma = \Deck(\tilde{X}/X)$ acts by positive homothethies of the \K \ metric. By \ref{thm:LCKacoperiri}, $X$ thus has an LCK form.

	\bigskip

	We conclude with the proof of the above claim. 
	
	\textbf{Proof of the Claim:} We say that a compact subset $L\subset X$ has property $(\star)$ if the following implication holds:
	$$L\cap \supp \chi_{(\alpha,\gamma)}\not = \emptyset \Rightarrow L\subset \supp\chi_{(\alpha,\gamma)}.$$
	For an $L$ with property $(\star)$, denote 
	$$index(L)=\{(\alpha,\gamma)\in A\times\Gamma: L\cap \supp\chi_{(\alpha,\gamma)}\not =\emptyset \},$$
	which is a nonempty, finite set. We fix an arbitrary point $x_0\in X$ and consider $L\subset X$ a compact neighborhood of $x_0$ with property $(\star)$. Then, $\tilde{g}(L)$ is also compact and $\tilde{g}(L)\subset\tilde{V}_{(\alpha,\gamma)}$ for all $(\alpha,\gamma)\in index(L)$. 
	
	Next, we want to prove the following assertion:
\begin{quote}[$\dagger$]\label{assertion}
there exists a constant $\delta=\delta_L$ such that $\phi_{(\alpha, \gamma)}^\epsilon$ is strongly plurisubharmonic on $L$ for any $0<\epsilon_\alpha\leq \delta$ with $\alpha\in index(L)$, and for any $\gamma\in\Gamma$. 
\end{quote}
	
	Now, if $L$ is chosen to be sufficiently small, then by using extensions for our functions in some local embeddings, we may consider, without loss of generality, that $\tilde{X}$ and $\tilde{Y}$ are open sets in the Euclidean spaces $\mathbb{C}^N$ and $\mathbb{C}^M$, respectively. 
	
	Further, for every $(\alpha,\gamma)\in index(L)$, there exist positive constants $p_\alpha$, $q_\alpha$, $b_\alpha$, $c_\alpha$ which depend only on $\alpha$, such that 
	\begin{enumerate}
	\item $\mathcal{L}(\psi_\alpha,\pi_Y(\tilde{y}))\xi\geq p_\alpha \|\xi\|^2$,
	\item $\mathcal{L}(\varphi_{(\alpha, \gamma)},\tilde{x})\zeta\geq \rho(\gamma) q_\alpha \|\zeta\|^2$,
	\item $|\varphi_{(\alpha, \gamma)}(\tilde{x})\mathcal{L}((\tau_\alpha\circ\pi_Y),\tilde{y})\xi|\leq \rho(\gamma) b_\alpha\|\xi\|^2$,
	\item $|\Re\partial(\tau_\alpha\circ\pi_Y)(\tilde{y})\xi \otimes \overline{\partial}\varphi_{(\alpha, \gamma)}(\tilde{x})\zeta|\leq \rho(\gamma) c_\alpha\|\xi\| \cdot\|\zeta\|$,
	\end{enumerate}
	where $\mathcal{L}$ denotes the Levi form, $\tilde{x}\in L$, $\tilde{y}=\tilde{g}(\tilde{x})$, $\zeta\in\mathbb{C}^N$, and $\xi=\partial \tilde{g}(\tilde{x})\zeta$. 
	(1) and (2) hold since $\psi_\alpha$ and $\varphi_{(\alpha, \gamma)}$ are strongly plurisubharmonic and $\tilde{g}(L)$, $L$ are compact. (3) and (4) are straightforward.
	
	Computing the Levi form of the function $\phi_{(\alpha, \gamma)}^\epsilon$ yields
	\begin{equation}\label{ineqLevi}
	\mathcal{L}(\phi_{(\alpha, \gamma)}^\epsilon,\tilde{x})\zeta\geq \rho(\gamma) (P\|\xi\|^2+ Q\|\zeta\|^2-B\|\xi\|^2-2C\|\xi\|\cdot\|\zeta\|),
	\end{equation}
	where we used the notations 
	\begin{itemize}
	\item $P=\min\{p_\alpha:\mbox{there exists }\gamma \mbox{ such that }(\alpha,\gamma)\in index(L)\},$
	\item $Q=Q(\tilde{x})=\sum \epsilon_\alpha q_\alpha \chi_{(\alpha,\gamma)}(\tilde{x})$,
	\item $B=B(\tilde{x})=2\sum \epsilon_\alpha b_\alpha \sqrt{\chi_{(\alpha,\gamma)}(\tilde{x})}$,
	\item $C=C(\tilde{x})=2\sum \epsilon_\alpha c_\alpha \sqrt{\chi_{(\alpha,\gamma)}(\tilde{x})}$.
	\end{itemize}
	
	If we assume for now that 
	\begin{equation}\label{condPQBC}
	(P-B)Q>C^2
	\end{equation}
	holds on $L$, then the right-hand side of $(\ref{ineqLevi})$ is positive for all $\zeta\not =0$. Hence, the existence of the previously mentioned $\delta=\delta_L$ is proved. 
	
	Now, in order to get (\ref{condPQBC}), choose $\delta_1>0$ small enough, such that for $0<\epsilon_\alpha\leq \delta_1$, we have $B\leq P/2$ on $L$. Then, if $0<\epsilon_\alpha\leq \delta_1$  for all $\alpha$ with $(\alpha,\gamma)\in index(L)$,  Schwartz's inequality leads to 
	$$(P-B)Q\geq \frac{P}{2}\sum \epsilon_\alpha q_\alpha \chi_{(\alpha,\gamma)}\geq P_1\left(\sum\sqrt{\epsilon_\alpha q_\alpha \chi_{(\alpha,\gamma)}}\right)^2.$$ 
	Since the inequality 
	$$\sqrt{\epsilon_\alpha q_\alpha}\geq \epsilon_\alpha\sqrt{\frac{q_\alpha}{\delta_1}}$$
	is always true, by choosing a convenient (possibly smaller) $\delta_1$, we obtain the desired inequality (\ref{condPQBC}).
	
	Returning to the main line of proof of our claim, fix $(\alpha,\gamma)\in A\times \Gamma$ and let $L_1,L_2,\ldots,L_m$ be  compact subsets of $X$ having the property $(\star)$ such that $(\alpha,\gamma)\in index(L_j)$ for all $j=1,2,\ldots,m$. Take $\delta_j>0$ according to assertion [$\dagger$] for $L=L_j$, $j=1,\ldots,m$, and set $\delta=\min\{ \delta_1,\ldots,\delta_m\}$.  Then, for any positive numbers $\{\epsilon_\alpha\}$ with $\epsilon_\alpha\leq\delta$, for $(\alpha,\gamma)\in index(L_j)$, $j=1,\ldots,m$ (and no other assumption on $\epsilon_\beta$ with $(\beta,\gamma)\not\in \cup_{j=1}^m index(L_j)$), the function $\phi_{(\alpha,\gamma)}^\epsilon$ is strongly plurisubharmonic on $\cup_{j=1}^{m}L_j$. Assuming we have taken $L_1,\ldots,L_m$ to cover $\supp \chi_{(\alpha,\gamma)}$, we obtain $\phi_{(\alpha,\gamma)}^\epsilon$ strongly plurisubharmonic on $\supp \chi_{(\alpha,\gamma)}$. Moreover, note that $\eta^*\phi_{(\alpha,\gamma)}^\epsilon=\rho(\eta)\phi_{(\alpha,\eta^{-1}\gamma)}^\epsilon$, hence $\phi_{(\alpha,\eta)}^\epsilon$ is also strongly plurisubharmonic on $\supp \chi_{(\alpha,\eta)}$ for all $\eta\in\Gamma$, without having to impose further conditions on $\{\epsilon_\alpha\}$. Finally, the local finiteness of $\{\supp \chi_{(\alpha,\gamma)}\}_{(\alpha,\gamma)\in A\times\Gamma}$ implies that only a finitely many conditions are imposed on every $\epsilon_\alpha$, thus concluding the proof of the \textbf{Claim}. 
\end{proof}

\bigskip

\textbf{Acknowledgments.} We would like to thank Liviu Ornea for his support and suggestions, Victor Vuletescu for helpful discussions about LCK spaces with singularities and Alexandra Otiman for carefully reading the paper and her comments.

\end{document}